\documentclass[runningheads,a4paper]{llncs}

\usepackage{amssymb,amsmath}   
\usepackage{comment}
\usepackage{IEEEtrantools}
\usepackage{tikz}
\newcommand\person[1]{\textsc{#1}}
\newcommand\Gaussian{\person{Gauss}ian}
\newcommand\Fourier{\person{Fourier}}
\newcommand\FourierSynthesis{\mathcal{F}}
\newcommand\FourierAnalysis{\mathcal{F}^{-1}}
\newcommand\DiscreteFourierAnalysis{\DiscreteFourierTransform^{-1}}
\newcommand\DiscreteFourierSynthesis{\DiscreteFourierTransform}

\newcommand\R{\mathbb{R}}
\newcommand\Q{\mathbb{Q}}
\newcommand\Z{\mathbb{Z}}
\newcommand\N{\mathbb{N}}
\DeclareMathOperator\DiscreteFourierTransform{DFT}
\DeclareMathOperator\dif{d}
\DeclareMathOperator\id{id}
\DeclareMathOperator\cis{cis}
\DeclareMathOperator\cisone{cis1}
\DeclareMathOperator\erf{erf}
\DeclareMathOperator\upper{upper}
\DeclareMathOperator\samplingrate{rate}

\newcommand \translater[2]{#2 \rightarrow #1}
\newcommand \dilate[2]{#2 \uparrow #1}
\newcommand \shrink[2]{#2 \downarrow #1}

\newcommand \abs[1]{\left\vert#1\right\vert}
\newcommand \norm[1]{\left\Vert#1\right\Vert}
\newcommand \eucnorm[1]{\norm{#1}_{2}}
\newcommand \lebesguenorm[2]{\norm{#2}_{#1}}
\newcommand \scalarprod[2]{\left\langle#1,#2\right\rangle}
\newcommand \lebesguespace[1]{\mathcal{L}_{#1}}
\newcommand \schwartzspace{\mathcal{S}}

\newcommand \conj[1]{\overline{#1}}
\newcommand \adjoint[1]{{#1}^{*}}
\newcommand \polyfunc[1]{\widehat{#1}}
\newcommand \charfunc[1]{\chi_{#1}}

\newcommand \keyword[1]{\emph{#1}}

\newcommand \eqnbreak{\noalign{\mbox{}}}
\newcommand \eqnlabel[1]{\yesnumber\label{eqn:#1}}    
\newcommand \seclabel[1]{\label{sec:#1}}
\newcommand \tablabel[1]{\label{tab:#1}}
\newcommand \figlabel[1]{\label{fig:#1}}

\newcommand \eqnref[1]{(\ref{eqn:#1})}          
\newcommand \secref[1]{Section~\ref{sec:#1}}    
\newcommand \tabref[1]{Table~\ref{tab:#1}}      
\newcommand \figref[1]{Figure~\ref{fig:#1}}     

\newcommand\figcaption[2]{\caption{\textit{#2}}\figlabel{#1}}
\newcommand\tabcaption[2]{\caption{\textit{#2}}\tablabel{#1}}

\usepackage{hyperref}
\urldef{\mailsa}\path|henning.thielemann@informatik.uni-halle.de|    
\newcommand{\keywords}[1]{\par\addvspace\baselineskip
\noindent\keywordname\enspace\ignorespaces#1}

\begin{document}

\mainmatter  

\title{An algebra for signal processing}


%
%
\author{Henning Thielemann}
%

\institute{Institut f\"{u}r Informatik,
Martin-Luther-Universit\"{a}t Halle-Wittenberg, Germany\\
\mailsa}

%
%

\maketitle

\begin{abstract}
Our paper presents an attempt to axiomatise signal processing.
Our long-term goal is to formulate signal processing algorithms
for an ideal world of exact computation
and prove properties about them,
then interpret these ideal formulations
and apply them without change to real world discrete data.
We give models of the axioms that are based on \Gaussian{} functions,
that allow for exact computations
and automated tests of signal algorithm properties.
\keywords{Algebra, Signal processing, Gaussian function}
\end{abstract}

\section{Introduction}
%

\subsection{Motivation}

In signal processing we consider real or complex valued functions.
These functions represent signals or frequency spectra
and their arguments are considered to be time values or frequency values,
respectively.
There are some fundamental operations
like pointwise multiplication ``$\cdot$'' and convolution ``$*$'' of signals.
The \Fourier{} transform $\FourierSynthesis$
converts between signals and frequency spectra.
For precise definitions of these operations see \secref{operations}.
Additionally there are some essential laws,
that every signal processing scientist is familiar with,
such as the law, that the \Fourier{} transform is an homomorphism,
that maps multiplication to convolution:
\[
\FourierSynthesis(x\cdot y) = \FourierSynthesis x * \FourierSynthesis y
.
\]
This is an important connection and it is amazingly simple,
but it is not quite true.

\begin{itemize}
\item
The first problem is,
that depending on the precise definitions of the involved operations
there may be a factor $2\pi$ or $\sqrt{2\pi}$ to make the above identity correct.
This is like working in the imperial system,
where a gallon is not just the cubic power of a length unit
but 231 cubic inch.
This is at least cumbersome and error-prone when done manually,
but it also complicates computer implementations.
With the transcendent factors we have to work in fields like
$\Q(\sqrt{2\pi})$ or $\Q(\pi)$
or would have to work with approximations,
but working with only rational numbers would of course
be still exact and more efficient.
However, if we succeed to suppress the transcendental factor
in the above identity by the change of a definition,
then it will show up in another place.

Thus our first contribution is
to select a set of operations that are common in signal processing
and list identities that show their interrelation in principle.
We use the degrees of freedom in the operation definitions
for simplifying the laws as much as possible,
that is, avoid constant factors and assert maximum symmetry
in \secref{operations}.
Metaphorically speaking, we try to define something
like a metric system for signal processing.

\item
The second problem is,
that the convolution and \Fourier{} transform are not always defined.
For instance the straightforward definitions
of the \Fourier{} transform and its inverse
as given in \secref{operations},
map from absolutely integrable functions to bounded functions,
i.e.
$\FourierSynthesis, \FourierAnalysis:
   \lebesguespace{1}(\R) \to \lebesguespace{\infty}(\R)
$
.
That is, in general you cannot apply the inverse \Fourier{} transform
to the result of a \Fourier{} transform.
There are even more such difficulties,
but they can be resolved by restriction
to the \person{Schwartz} space~$\schwartzspace$
of arbitrarily smooth and rapidly decaying functions.

\item
The third problem is,
that commonly signal processing in a computer is performed
on discretised functions with finite precision numbers.
That is, laws as the above one do not hold exactly.

Our contribution concerning this problem
is in the subsections of \secref{gaussians}.
There we take the presented laws as axioms
and develop models in terms of extended \Gaussian{} functions.
All of these functions are in $\schwartzspace$.
The more we extend those functions by factors and parameters,
the more operations we can support.
Our goal is to give definitions of functions,
such that we only need to cope with rational parameters.
This way we can represent a wide range of functions,
we can compute exactly
and we could even generalise these models to any algebraic field.




\end{itemize}

\subsection{Basics}

The theory of Algebra provides
many axiomatically described structures,
that generalise common mathematical objects and operations on them.
E.g. groups abstract permutations,
lattices abstract \person{Boole}an logic,
rings abstract integer arithmetic,
fields abstract rational arithmetic.
If we are able to perform a proof
with the axioms of a particular algebraic structure,
then our proof automatically
applies to every such structure.
In this paper we would like to abstract
from what is commonly called signal processing.

Signal processing is the research area
of construction, transformation and analysis
of oscillation functions in one variable (the time),
approximation of those functions and related operations with discretised data
and efficient implementations on digital computers.
An important view on an audio signal is the frequency spectrum,
because that is closer to the way humans hear
than the time-domain representation of a signal.
We obtain a frequency spectrum by the \Fourier{} transform of a signal.
Closely related to signal processing is image processing,
with the main difference being,
that functions in two variables are studied.
Another related area are random distributions in stochastics.
The convolution, a common operation in signal processing,
of two random distributions yields the distribution
of the sum of two random variables.
The \Fourier{} transform yields the characteristic function
of a random distribution
and there are more connections (see \secref{gaussians}).

However, signal processing, image processing and random distributions
are relatively seldom viewed from an algebraic point of view.
This may have to do with the tradition to focus on signal values (e.g. $f(t)$)
rather than on larger objects like signals (e.g. $f$).
E.g. it is common to express the delay of a signal by $f(t-d)$
rather than to use an operator like in $\translater{d}{f}$.
On the one hand this has the advantage, that some properties are obvious
(e.g. $f(t-(d_0+d_1)) = f((t-d_0)-d_1)$),
since they are a consequence of simple arithmetic,
but are not obvious at the higher level
(i.e. $\translater{(d_0+d_1)}{f} = \translater{d_0}{(\translater{d_1}{f})}$).
On the other hand we cannot cleanly express identities
involving intrinsically functional transforms
like the \Fourier{} transform,
that has lead to custom notations like
$F(\omega) \stackrel{\FourierSynthesis}{\leftrightarrow} f(t)$,
that could be expressed the functional way as
$\FourierSynthesis{F} = f$.

\section{Finding a set of fundamental operations}
\seclabel{operations}

We like to start this section listing the operations,
that are commonly used in signal processing,
together with typical applications.
We begin with the operations with indisputable definitions
and continue with the ones,
where differing definitions are around in the literature.
For the variant definitions we show, what laws they imply.
Then we choose the definitions
that make the laws most simple.
We close the section with a comprehensive list of laws
that hold for our definitions.

In \tabref{operations} we list the signal operations
together with their definitions
and in \figref{graphs} we illustrate them using example signals.
\begin{itemize}
\item Shrinking a signal as in \eqnref{shrink-graph}
means to increase all contained frequencies proportionally
and to shorten time accordingly.
\item Translating a signal as in \eqnref{translation-graph}
means to delay it.
\item Summing two signals as in \eqnref{sum-graph}
means to superpose them, that is, to play them simultaneously.
\item Multiplying two signals as in \eqnref{multiplication-graph}
means to control the amplitude of one signal by the other one.
For certain choices of signals this is also known as ring-modulation.
\item Convolving two signals as in \eqnref{convolution-graph}
means to apply the sound characteristics of one signal to the other one.
It may be used for suppressing or emphasising certain frequencies,
for smoothing or for application of reverb.
\end{itemize}

\begin{table}
\tabcaption{operations}{Basic signal processing operations}
\begin{IEEEeqnarray*}{s's'rCl}
operation & application & \mbox{definition} \\
\hline\\
shrinking & alter pitch and time &
   (\shrink{k}{x})(t) &=& x(k\cdot t) \\
translate & delay &
   (\translater{k}{x})(t) &=& x(t-k) \\
adjoint & &
   \adjoint{x} &=& \shrink{-1}{\conj{x}} \\
sum & mixing &
   (x+y)(t) &=& x(t) + y(t) \\
multiplication & envelope &
   (x\cdot y)(t) &=& x(t) \cdot y(t) \\
convolution & frequency filter &
   (x*y)(t) &=& \int_{\R} x(\tau)\cdot y(t-\tau) \dif \tau \\
modulation & frequency shift &
   x \cdot \cisone & & \mbox{where } \cisone t = \exp(2\pi i \cdot t)
\end{IEEEeqnarray*}
\end{table}
%
%
\newenvironment{ctikzpicture}{
\begin{tabular}[c]{c}
\begin{tikzpicture}
}
{
\end{tikzpicture}
\end{tabular}
}

\begin{figure}
\begin{IEEEeqnarray*}{rCl}
\shrink{2}{
\begin{ctikzpicture}
\draw[->] (-1.6,0) -- (1.6,0);
\draw[->] (0,-0.2) -- (0,1.2);
\draw (-0.5,0) -- (-0.5,1) -- (0.5,1) -- (0.5,0);
\end{ctikzpicture}
}
&=&
\begin{ctikzpicture}
\draw[->] (-1.6,0) -- (1.6,0);
\draw[->] (0,-0.2) -- (0,1.2);
\draw (-0.25,0) -- (-0.25,1) -- (0.25,1) -- (0.25,0);
\end{ctikzpicture}
\eqnlabel{shrink-graph}
\\
\translater{\frac{1}{2}}{
\begin{ctikzpicture}
\draw[->] (-1.6,0) -- (1.6,0);
\draw[->] (0,-0.2) -- (0,1.2);
\draw (-0.8,0) -- (-0.8,1) -- (0.2,1) -- (0.2,0);
\end{ctikzpicture}
}
&=&
\begin{ctikzpicture}
\draw[->] (-1.6,0) -- (1.6,0);
\draw[->] (0,-0.2) -- (0,1.2);
\draw (-0.3,0) -- (-0.3,1) -- (0.7,1) -- (0.7,0);
\end{ctikzpicture}
\eqnlabel{translation-graph}
\\
\begin{ctikzpicture}
\draw[->] (-1.6,0) -- (1.6,0);
\draw[->] (0,-0.7) -- (0,0.7);
\draw (-1.0,0) -- (-1.0,-0.5) -- (0.5,-0.5) -- (0.5,0);
\end{ctikzpicture}
+
\begin{ctikzpicture}
\draw[->] (-1.6,0) -- (1.6,0);
\draw[->] (0,-0.7) -- (0,0.7);
\draw (-0.5,0) -- (-0.5,0.5) -- (1.0,0.5) -- (1.0,0);
\end{ctikzpicture}
&=&
\begin{ctikzpicture}
\draw[->] (-1.6,0) -- (1.6,0);
\draw[->] (0,-0.7) -- (0,0.7);
\draw (-1.0,0) -- (-1.0,-0.5) -- (-0.5,-0.5) -- (-0.5,0)
   -- ( 0.5,0) -- ( 0.5, 0.5) -- ( 1.0, 0.5) -- ( 1.0,0);
\end{ctikzpicture}
\eqnlabel{sum-graph}
\\
\begin{ctikzpicture}
\draw[->] (-1.6,0) -- (1.6,0);
\draw[->] (0,-0.7) -- (0,1.2);
\draw (-1.5,0) -- (-1.0, 0.5) -- (-1.0,-0.5) -- (0.0, 0.5)
               -- ( 0.0,-0.5) -- ( 1.0, 0.5) -- (1.0,-0.5) -- (1.5,0);
\end{ctikzpicture}
\cdot
\begin{ctikzpicture}
\draw[->] (-1.6,0) -- (1.6,0);
\draw[->] (0,-0.7) -- (0,1.2);
\draw (-0.5,0) -- (-0.5,1) -- (0.5,1) -- (0.5,0);
\end{ctikzpicture}
&=&
\begin{ctikzpicture}
\draw[->] (-1.6,0) -- (1.6,0);
\draw[->] (0,-0.7) -- (0,1.2);
\draw (-0.5,0) -- (0,0.5) -- (0,-0.5) -- (0.5,0);
\end{ctikzpicture}
\eqnlabel{multiplication-graph}
\\
\begin{ctikzpicture}
\draw[->] (-1.6,0) -- (1.6,0); 
\draw[->] (0,-0.2) -- (0,1.2); 
\draw (-0.5,0) -- (-0.5,1) -- (0.5,1) -- (0.5,0);
\end{ctikzpicture}
*
\begin{ctikzpicture}
\draw[->] (-1.6,0) -- (1.6,0); 
\draw[->] (0,-0.2) -- (0,1.2); 
\draw (-0.5,0) -- (-0.5,1) -- (0.5,1) -- (0.5,0);
\end{ctikzpicture}
&=&
\begin{ctikzpicture}
\draw[->] (-1.6,0) -- (1.6,0); 
\draw[->] (0,-0.2) -- (0,1.2); 
\draw (-1,0) -- (0,1) -- (1,0);
\end{ctikzpicture}
\eqnlabel{convolution-graph}
\\
\FourierSynthesis~
\begin{ctikzpicture}
\draw[->] (-1.6,0) -- (1.6,0);
\draw[->] (0,-0.4) -- (0,1.2);
\draw (-1,0) -- (-1,0.5) -- (1,0.5) -- (1,0);
\end{ctikzpicture}
&=&
\begin{ctikzpicture}
\draw[->] (-1.6,0) -- (1.6,0);
\draw[->] (0,-0.4) -- (0,1.2);
\draw plot file {sinc.csv};
\end{ctikzpicture}
\eqnlabel{fourier-graph}
\\
\end{IEEEeqnarray*}
\figcaption{graphs}{
Illustration of the basic signal processing operations.
}
\end{figure}

We continue with the controversial definitions.
Actually, there is essentially one operation,
that is defined differently throughout the signal processing literature:
The \Fourier{} transform.
It is the transform that computes the frequency spectrum
for a signal in the time domain (\keyword{backward} or \keyword{analysis} transform)
and vice versa (\keyword{forward} or \keyword{synthesis} transform).
For an example see \eqnref{fourier-graph} in \figref{graphs}.
The other definition with varying instances
is the one for functions that are eigenfunctions of the \Fourier{} transform.
Obviously it depends entirely on the definition of the \Fourier{} transform.
\begin{enumerate}
\item Oscillations with period $1$
  \begin{IEEEeqnarray*}{s'rCl}
  \Fourier{} forward &
    \FourierSynthesis_1 x(\tau)
       &=& \int_\R \exp(2\pi i\cdot \tau\cdot t)\cdot x(t) \dif t \\
  \Fourier{} backward &
    \FourierAnalysis_1 x(\tau)
       &=& \int_\R \exp(-2\pi i\cdot \tau\cdot t)\cdot x(t) \dif t \\
  duality &
    \FourierSynthesis_1 (\FourierSynthesis_1 x) &=& \shrink{-1}{x} \\
  eigenfunction &
    g(t) &=& \exp\left(-\pi\cdot t^2\right) \eqnlabel{eigen1} \\
  eigenvalue is 1 &
    \FourierSynthesis_1 g &=& g \\
  unitarity &
    \scalarprod{x}{y}
       &=& \scalarprod{\FourierSynthesis_1 x}{\FourierSynthesis_1 y} \\
  convolution theorem &
    \FourierSynthesis_1(x*y)
       &=& \FourierSynthesis_1 x \cdot \FourierSynthesis_1 y \\
  &
    \FourierSynthesis_1(x\cdot y)
       &=& \FourierSynthesis_1 x * \FourierSynthesis_1 y \\
  derivative &
    %
    \FourierSynthesis_1(x')
       &=& -2\pi i\cdot\id \cdot \FourierSynthesis_1 x
              \eqnlabel{derivative1}
  \end{IEEEeqnarray*}
\item Oscillations with period $2\pi$
  \begin{IEEEeqnarray*}{s'rCl}
  \Fourier{} forward &
    \FourierSynthesis_2 x(\tau)
       &=& \frac{1}{\sqrt{2\pi}}\cdot
               \int_\R \exp(i\cdot \tau\cdot t)\cdot x(t) \dif t \\
  \Fourier{} backward &
    \FourierAnalysis_2 x(\tau)
       &=& \frac{1}{\sqrt{2\pi}}\cdot
               \int_\R \exp(-i\cdot \tau\cdot t)\cdot x(t) \dif t \\
  duality &
    \FourierSynthesis_2 (\FourierSynthesis_2 x) &=& \shrink{-1}{x} \\
  eigenfunction &
    g(t) &=& \exp\left(-t^2\right) \eqnlabel{eigen2} \\
  eigenvalue is 1 &
    \FourierSynthesis_2 g &=& g \\
  unitarity &
    \scalarprod{x}{y}
       &=& \scalarprod{\FourierSynthesis_2 x}{\FourierSynthesis_2 y} \\
  convolution theorem &
    \FourierSynthesis_2(x*y)
       &=& \sqrt{2\pi}\cdot \FourierSynthesis_2 x \cdot \FourierSynthesis_2 y
              \eqnlabel{convolution2} \\
  &
    \FourierSynthesis_2(x\cdot y)
       &=& \frac{1}{\sqrt{2\pi}}\cdot \FourierSynthesis_2 x * \FourierSynthesis_2 y
              \eqnlabel{multiplication2} \\
  derivative &
    \FourierSynthesis_2(x')
       &=& -i\cdot\id \cdot \FourierSynthesis_2 x
              \eqnlabel{derivative2}
  \end{IEEEeqnarray*}
\item Oscillations with period $2\pi$ and no roots of $\pi$
  \begin{IEEEeqnarray*}{s'rCl}
  \Fourier{} forward &
    \FourierSynthesis_3 x(\tau)
       &=& \int_\R \exp(i\cdot \tau\cdot t)\cdot x(t) \dif t \\
  \Fourier{} backward &
    \FourierAnalysis_3 x(\tau)
       &=& \frac{1}{2\pi}\cdot
             \int_\R \exp(-i\cdot \tau\cdot t)\cdot x(t) \dif t \\
  \end{IEEEeqnarray*}
\end{enumerate}

We do not further consider Definition 3,
because forward and backward \Fourier{} transform have different factors,
thus laws for forward and backward transform differ in factors.
For definitions 1 and 2 many laws are equal up to the choice of the transform direction.

At the first glance the definitions 1 and 2 of the \Fourier{} transform
seem to be equally convenient or equally inconvenient.
Thus many textbooks just introduce a definition
and do not explain, why they prefer the one to the other possible ones.
However, we think that factors in laws are bad,
especially worse than factors that can be hidden in a function definition.
E.g.\ compare the eigenfunction property of the \Gaussian{} function
in \eqnref{eigen1} and \eqnref{eigen2}:
We can define the \Gaussian{} bell curve in both ways and call it just $g$.
The factor $\pi$ does no longer get in the way
when transforming equations containing $\FourierSynthesis$ and $g$.
The same applies to the laws on derivatives
in \eqnref{derivative1} and \eqnref{derivative2},
where we can consider $2\pi i\cdot\id$ as one function.
In contrast to that, the convolution theorems
in \eqnref{convolution2} and \eqnref{multiplication2}
have constant factors, even different ones.
This makes manual equation manipulation cumbersome and error-prone.
Even more it makes computations exclusively with rational numbers impossible,
due to the irrational factor $\sqrt{2\pi}$.
We could suppress the factors
in \eqnref{convolution2} and \eqnref{multiplication2}
by defining convolution or multiplication containing a constant factor,
but we think this is not natural.

For these reasons we will stick to definition 1
and omit the index of $\FourierSynthesis$ in the rest of this paper.
%
%
\seclabel{axioms}
\newcommand\redundant{\mbox{ }^\#}
Below we list the laws that follow from these definitions
and that we want as axioms for our algebraic structure.
Be aware, that we have omitted convergence conditions
for the operations that involve integration.
The laws marked with \# 
can be derived from the remaining laws
and are just given for the purpose of completeness.
\begin{IEEEeqnarray*}{rCl"rCl}
x+y &=& y+x &
x+(y+z) &=& (x+y)+z \\
x\cdot y &=& y\cdot x &
x\cdot (y \cdot z) &=& (x\cdot y) \cdot z \\ &&&
x\cdot (y+z) &=& x\cdot y + x \cdot z \\
x*y &=& y*x\redundant &
x*(y*z) &=& (x*y)*z\redundant \\ &&&
x*(y+z) &=& x*y + x*z\redundant \\
\translater{0}{x} &=& x &
\translater{b}{(\translater{a}{x})} &=& \translater{(a+b)}{x} \\
\shrink{1}{x} &=& x &
\shrink{b}{(\shrink{a}{x})} &=& \shrink{(a\cdot b)}{x} \\ &&&
\translater{a}{(\shrink{b}{x})} &=& \shrink{b}{(\translater{(a\cdot b)}{x})} \\
\translater{a}{(x+y)} &=& (\translater{a}{x}) + (\translater{a}{y}) &
\shrink{a}{(x+y)} &=& (\shrink{a}{x}) + (\shrink{a}{y}) \\
\translater{a}{(x\cdot y)} &=& (\translater{a}{x}) \cdot (\translater{a}{y}) &
\shrink{a}{(x\cdot y)} &=& (\shrink{a}{x}) \cdot (\shrink{a}{y}) \\
\translater{a}{(x*y)} &=& x * (\translater{a}{y}) &
\shrink{a}{(x*y)} &=& \abs{a} \cdot (\shrink{a}{x}) * (\shrink{a}{y}) \\
(x\cdot y)' &=& x'\cdot y + x \cdot y' &
(x * y)' &=& x * y' \\
\eqnbreak
\FourierSynthesis(x+y) &=& \FourierSynthesis x + \FourierSynthesis y &
\FourierSynthesis(k\cdot x) &=& k \cdot \FourierSynthesis x \\
\FourierSynthesis(x*y) &=& \FourierSynthesis x \cdot \FourierSynthesis y &
\FourierSynthesis(x\cdot y) &=& \FourierSynthesis x * \FourierSynthesis y\redundant \\
\scalarprod{x}{y}
         &=& \scalarprod{\FourierSynthesis x}{\FourierSynthesis y} &
\eucnorm{x} &=& \eucnorm{\FourierSynthesis x}\redundant \\
\FourierSynthesis(\shrink{a}{x})
         &=& \abs{a}\cdot\shrink{\frac{1}{a}}{\FourierSynthesis x} &
\FourierSynthesis(\translater{a}{x})
         &=& (\shrink{a}{\cis}) \cdot \FourierSynthesis x \\
\FourierSynthesis(\FourierSynthesis x)
         &=& \shrink{-1}{x} \\
\FourierSynthesis(\adjoint{x})
         &=& \conj{\FourierSynthesis{x}} &
\FourierSynthesis (x') &=& - 2\pi i \cdot \id \cdot \FourierSynthesis x
\end{IEEEeqnarray*}

\section{Development of Gaussian models}
\seclabel{gaussians}

Now that we have stated some axioms,
we want to construct an ideal world,
that is, a class of functions (or signals),
where they hold.
This class of function shall allow for simple constraints of the laws,
for exact and efficient computation
(more precisely: operations in a field),
and shall allow to represent many mathematically important objects exactly
and real world signals approximately.

Since \Gaussian{} functions are eigenfunctions of the \Fourier{} transform,
they are perfect for representing signals in both time and frequency domain.
They can be extended in a relatively simple way,
such that all of the operations can be performed, that we listed initially.

We like to stress that \Gaussian{} functions are not only interesting,
because they allow to perform the operations we want,
but \Gaussian{} functions are central to signal processing and stochastics.
\begin{itemize}
\item \Gaussian{} functions are used as filter window for smoothing.
\item \Gaussian{} functions ``minimise uncertainty'',
that is, they are the functions
where \person{Heisenberg}'s uncertainty relation
becomes an equation.
\cite{triebel1980analysis}
\item With complex modulation \Gaussian{} functions
are called \person{Gabor} atoms
or \keyword{time-frequency atoms} in the \person{Gabor} transform.
The \person{Gabor} transform is a windowed \Fourier{} transform,
that shows how frequency components evolve over time.
\item Complex modulated \Gaussian{} functions
with a correction offset are called \person{Morlet} wavelets
and are used in the Continuous Wavelet Transform.
This transform is also intended for showing the evolution
of frequency components over time,
but it has higher time resolution and less frequency resolution
for high frequencies.
\item Best basis pursuits and matching pursuits are techniques
for decomposing a signal
into a finite number of irregularly located time-frequency atoms.
\cite{mallat1993matchingpursuits}
\end{itemize}
Since \person{Gabor} transform, \person{Morlet} wavelet transform
as well as best basis and matching pursuits
aim at decomposition of a signal into time-frequency atoms,
we have several tools for approximating real world signals using
those atoms as building blocks.

In stochastics the density of the Normal distribution is a \Gaussian{} function.
The Central Limit Theorem states,
that adding more and more random variables
(and divide by the square root of the number of added variables)
in most practical cases approaches a normally distributed random variable.
Translated to signal processing this means,
that smoothing a signal again and again with practically any filter window,
approximates a \Gaussian{} filter.

The derivation of the representations below is not particularly difficult,
but we focus on finding representations that simplify most operations
in terms of use of transcendental constants and irrational algebraic functions.
We have implemented these function classes in Haskell
using the NumericPrelude type class hierarchy.
You find our implementation at \\
\centerline{
\url{http://code.haskell.org/numeric-prelude/src/MathObj/Gaussian/}.
}

\subsection{Simple \Gaussian{}s}
\seclabel{gaussians1}

Simple \Gaussian{} functions shall be the first class of functions
that we want to consider.
In order to avoid a transcendental factor containing $\pi$
in any of our function parameters
when applying \Fourier{} transform $\FourierSynthesis$,
we cannot choose $f(t)=\exp(-t^2)$
but we have to choose its eigenfunction
$f(t)=\exp(-\pi\cdot t^2)$.
We also want to support shrinking of a function,
what requires adding a shrinking parameter~$c$,
that we like to write in \keyword{curried form}:
$f(c)(t)=\exp(-\pi\cdot (c\cdot t)^2)$.
However this would yield a square root
in the convolution of two such functions.
It can be prevented by choosing the form
$f(c)(t)=\exp(-\pi\cdot c\cdot t^2)$ with $c\in\Q$.
The \Fourier{} transform of a shrunk function
leads to another factor,
the amplitude $y$ of the function.
We end up with the function form:
\begin{IEEEeqnarray*}{rCl}
f(y,c)(t) &=& \sqrt{y}\cdot\exp(-\pi\cdot c\cdot t^2) \\
c \in \Q && y \in [0,\infty) \cap \Q
.
\end{IEEEeqnarray*}

This simple function class already allows for several operations,
where convolution and \Fourier{} transform have constraints
that assert convergence:
\begin{IEEEeqnarray*}{s'rCl}
scaling &
   k \cdot f(y,c) &=& f(y\cdot k^2, c) \\
shrinking &
   \shrink{k}{f(y,c)} &=& f(y, c \cdot k^2) \\
conjugate &
   \conj{f(y,c)} &=& f(y, c) \\
multiplication &
   f(y_0,c_0) \cdot f(y_1,c_1) &=& f(y_0\cdot y_1, c_0+c_1) \\
power with $r\ge 0$ &
   f(y,c)^r &=& f(y^r, r\cdot c) \\
convolution $c_0 + c_1 > 0$ &
   f(y_0,c_0) * f(y_1,c_1) &=&
      f\left(\frac{y_0\cdot y_1}{c_0 + c_1}, \frac{c_0\cdot c_1}{c_0 + c_1}\right) \\
\Fourier{} transform $c>0$ &
   \FourierAnalysis{(f(y,c))} &=&
      f\left(\frac{y}{c}, \frac{1}{c}\right)
.
\end{IEEEeqnarray*}

Typical functionals like function norms do not easily satisfy our goal
of using the most simple algebraic structures.
They need roots of parameters
or constant transcendental factors,
and thus require special treatment:
\begin{IEEEeqnarray*}{s'rCl}
$\lebesguespace{1}$-norm &
   \norm{f(y,c)}_1 &=& \sqrt{\frac{y}{c}} \\
$\lebesguespace{2}$-norm &
   \norm{f(y,c)}_2 &=& \sqrt{\frac{y}{\sqrt{2\cdot c}}} \\
$\lebesguespace{\infty}$-norm &
   \norm{f(y,c)}_{\infty} &=& \sqrt{y} \\
$\lebesguespace{p}$-norm &
   \norm{f(y,c)}_{p} &=& \sqrt{\frac{y}{\sqrt[p]{p\cdot c}}} \\
variance &
   \frac
     {\lebesguenorm{1}{t \mapsto t^2\cdot f(y,c)(t)}}
     {\lebesguenorm{1}{f(y,c)}}
    &=& \frac{1}{2\pi\cdot c}
.
\end{IEEEeqnarray*}

\subsection{Translated and modulated \Gaussian{}s}
\seclabel{gaussians2}

In the next step we want to translate and modulate the \Gaussian{} function.
The most simple function class, that allows this, seems to be:
\begin{IEEEeqnarray*}{rCl}
f(y,a,b,c)(t) &=&
 \sqrt{y}\cdot\exp\left(-\pi\cdot (a+b\cdot t+c\cdot t^2)\right)
   \eqnlabel{gaussian-translated} \\
&& y \in [0,\infty) \cap \Q \\
c \in \Q && \{a,b\} \subset \Q+i\Q
.
\end{IEEEeqnarray*}
With this representation we can perform the following operations:
%
%
\newcommand\opname[1]{\noalign{#1}}%
\newcommand\opnameinline[1]{\mbox{#1\quad}\hfill}%
\begin{IEEEeqnarray*}{rCl}
\opnameinline{translation}
   \translater{k}{f(y,a,b,c)} &=&
      f(y, a - b\cdot k + c\cdot k^2,
           b - 2\cdot c\cdot k, c) \\
\opnameinline{modulation}
   f(y,a,b,c) \cdot (\shrink{k}{\cisone}) &=&
      f\left(y, a, b+2\cdot i\cdot k, c\right) \\
\opnameinline{scaling}
   k \cdot f(y,a,b,c) &=& f(y\cdot k^2, a,b,c) \\
\opnameinline{shrinking}
   \shrink{k}{f(y,a,b,c)} &=& f(y, a, b\cdot k, c \cdot k^2) \\
\opnameinline{conjugate}
   \conj{f(y,a,b,c)} &=& f(y, \conj{a}, \conj{b}, c) \\
\opnameinline{power with $n\ge 0$}
   f(y,a,b,c)^n &=& f(y^n, n\cdot a, n\cdot b, n\cdot c) \\
\opname{multiplication}
   f(y_0,a_0,b_0,c_0) \cdot f(y_1,a_1,b_1,c_1)
      &=& f(y_0\cdot y_1, a_0+a_1, b_0+b_1, c_0+c_1) \\
\opname{convolution $c_0 + c_1 > 0$}
   f(y_0,a_0,b_0,c_0) * f(y_1,a_1,b_1,c_1) &=&
 \\ \noalign{\hfill
      $f\left(\frac{y_0\cdot y_1}{c_0 + c_1},
             a_0 + a_1 - \frac{(b_0 - b_1)^2}{4\cdot(c_0 + c_1)},
             \frac{b_0\cdot c_1 + b_1\cdot c_0}{c_0 + c_1},
             \frac{c_0\cdot c_1}{c_0 + c_1}\right)$ }
\opname{\Fourier{} transform $c>0$}
   \FourierAnalysis{(f(y,a,b,c))} &=&
      f\left(\frac{y}{c}, a-\frac{b^2}{4c}, -\frac{i\cdot b}{c}, \frac{1}{c}\right)
\end{IEEEeqnarray*}
The correctness of the equations for the \Fourier{} transform
and the convolution are not so obvious.
The \Fourier{} transform can be derived by translating the function to the origin
and demodulate it, such that it becomes real.
Then do \Fourier{} transform and
translate and modulate it corresponding to the normalisations
that we performed in time domain.
The convolution can also be derived from such an normalisation.
An alternative is to multiply in frequency domain.

We could also employ the definition
\begin{IEEEeqnarray*}{rCl}
f(y,a,b,c)(t) =
    \sqrt{y}\cdot\exp(- (a+b\cdot\sqrt\pi t+c\cdot \pi t^2))
   \eqnlabel{gaussian-translated-sqrt-pi}
\end{IEEEeqnarray*}
and the formulas for most operations would remain the same,
but translation, shrinking and modulation would have to be interpreted
with respect to the unit~$\sqrt\pi$.

With the considered representation
we can also represent two other kinds of functions
that are important to signal processing:
A decaying exponential curve can be obtained with $c=0 \land b>0$.
It is frequently encountered as envelope of percussive sounds.
Unfortunately that curve is unrestricted in time,
whereas in natural sounds the envelope usually starts suddenly somewhere in time.

The other important signal is a tone of linearly changing frequency,
a \keyword{chirp}.
In order to represent it,
we have to generalise the real parameter $c$ to a complex parameter.
Setting $\Im(c) \ne 0, \Re(c) = 0$ yields a pure chirp,
whereas $\Re(c) > 0$ yields a \keyword{chirplet},
that is a chirp, that is localised in time.
The name is chosen analogously to \keyword{wavelets}.
The \Fourier{} transforms and norm computations converge,
if and only if, $\Re{c}>0$ (i.e. only for chirplets),
and the convolution exists if and only if $\Re(c_0+c_1)>0$.
The \person{Bluestein} transform allows us
to write a \Fourier{} transform in terms of a convolution with a chirp.
\[
\FourierAnalysis{x} = ((x\cdot f(1,0,0,-i)) * f(1,0,0,i)) \cdot f(1,0,0,-i)
\]
The downside of this generalisation is,
that we need to maintain a complex amplification factor $y$.
This involves a complex square root
and we must choose the right branch.
For a single convolution or \Fourier{} transform
choosing the branch with positive real part is just the right thing.
The real part cannot vanish, since then the integrals do not converge at all.
But when multiplying square roots
with the convention of positive real parts
then we must respect
\begin{IEEEeqnarray*}{rCl}
\sqrt{c_0} \cdot \sqrt{c_1} &=& (-1)^k \cdot \sqrt{c_0\cdot c_1} \\
k &=&
\begin{cases}
%
1 : (\upper{c_0} \land \upper{c_1} \land \neg \upper(c_0\cdot c_1)) \lor \\
\quad (\neg \upper{c_0} \land \neg \upper{c_1} \land \upper(c_0\cdot c_1)) \\
0 : \mbox{otherwise}
\end{cases}
\\
\upper c &=& \Im{c} > 0 \lor (\Im{c} = 0 \land \Re{c} < 0)
.
\end{IEEEeqnarray*}
That is we must maintain the sign of the real part separately.
In the representation \eqnref{gaussian-translated}
we can implement a flipped sign by adding $i$ to $b$.
However maintaining this sign means comparisons
and thus would not work for generalisations from $\Q$ to other fields,
e.g. finite fields.

\subsection{\Gaussian{}s multiplied with polynomials}
\seclabel{gaussians3}

Another important operation in signal processing is derivation,
be it in differential equations like oscillation equations,
for the representation of an eigenbasis of the \Fourier{} transform,
or as a \keyword{highpass filter},
that is, a frequency filter that emphasises high frequencies
and suppresses low frequencies.
Derivation requires, that we extend our representation
to a product of a \Gaussian{} function and a polynomial function.
However using the simple representation
$f(y,a,b,c)(t) =
   \sqrt{y}\cdot\exp\left(-\pi\cdot (a+b\cdot t+c\cdot t^2)\right)$
from \eqnref{gaussian-translated},
this would mean to maintain polynomial expressions of $\pi$
as coefficients of the polynomial factor.

We like to avoid that and instead extend \eqnref{gaussian-translated-sqrt-pi} to:
\begin{IEEEeqnarray*}{rCl}
f(y,a,b,c)(t) &=&
    \sqrt{y}\cdot\exp\left(- (a+b\sqrt\pi\cdot t+c\cdot \pi\cdot t^2)\right) \\
\varphi((y,a,b,c),p)(t) &=&
    f(y,a,b,c)(t) \cdot \polyfunc{p}(t) \\
&& y \in [0,\infty) \cap \Q \\
c \in \Q && \{a,b\} \subset \Q+i\Q \\
\polyfunc{p}(t) &\in& (\Q+i\Q)[\sqrt{\pi} \cdot t] \\
\polyfunc{p}(t) &=& \sum_{j=0}^{n} p_j \cdot (\sqrt{\pi} \cdot t)^j
\end{IEEEeqnarray*}
As mentioned above translation and modulation have to be interpreted
with respect to a unit $\sqrt\pi$,
and derivation must contain a factor of $\sqrt\pi$.
However in order to get the usual units
we can still replace $\Q$ by $\Q(\sqrt{\pi})$.

For the following list of instantiations of signal processing transforms
we like to subsume the parameters of $f$ in a parameter tuple $\alpha$.
\begin{IEEEeqnarray*}{rCl}
\opnameinline{translation}
   \translater{\sqrt\pi k}{\varphi(\alpha,p)} &=&
      (\translater{\sqrt\pi k}{f\alpha}) \cdot
      (\translater{\sqrt\pi k}{\polyfunc{p}}) \\
\opnameinline{modulation}
   \varphi(\alpha,p) \cdot \left(\shrink{\frac{k}{\sqrt\pi}}{\cisone}\right) &=&
      \left(f\alpha \cdot \left(\shrink{\frac{k}{\sqrt\pi}}{\cisone}\right)\right) \cdot
      \polyfunc{p} \\
\opnameinline{scaling}
   k \cdot \varphi(\alpha,p) &=&
      f\alpha \cdot (k\cdot\polyfunc{p}) \\
\opnameinline{shrinking}
   \shrink{k}{\varphi(\alpha,p)} &=&
      (\shrink{k}{f\alpha}) \cdot
      (\shrink{k}{\polyfunc{p}}) \\
\opnameinline{conjugate}
   \conj{\varphi(\alpha,p)} &=&
      \conj{f\alpha} \cdot \polyfunc{\conj{p}} \\
\opnameinline{multiplication}
   \varphi(\alpha_0,p_0) \cdot \varphi(\alpha_1,p_1) &=&
      (f(\alpha_0) \cdot f(\alpha_1)) \cdot
      \polyfunc{p_0 \cdot p_1} \\
\opnameinline{power with $n\in\N$}
   \varphi(\alpha,p)^n &=&
      (f\alpha)^n \cdot \polyfunc{p^n} \\
\opnameinline{convolution}
   \varphi(\alpha_0,p_0) * \varphi(\alpha_1,p_1) &=&
      \FourierSynthesis\left(
         \FourierAnalysis{(\varphi(\alpha_0,p_0))}\cdot
         \FourierAnalysis{(\varphi(\alpha_1,p_1))}
      \right)
\end{IEEEeqnarray*}
\begin{IEEEeqnarray*}{rCl}
\opname{\Fourier{} transform}
   \FourierAnalysis{(\varphi(\alpha,s:p))} &=&
      s\cdot\FourierAnalysis{(f\alpha)} +
      \frac{i}{2\cdot\sqrt\pi}\cdot
        \left(\FourierAnalysis{(\varphi(\alpha,p))}\right)' \\
   && \qquad \mbox{where}\quad
         \polyfunc{s:p}(t) = s + \sqrt\pi\cdot t\cdot\polyfunc{p}(t) \\
\opname{Differentiation}
   \frac{1}{\sqrt\pi}\cdot(\varphi((y,a,b,c),p))' &=&
      f(y,a,b,c)\cdot \left(
         \frac{1}{\sqrt\pi}\cdot\polyfunc{p}' -
         (t \mapsto b + c\cdot \sqrt\pi\cdot t) \cdot \polyfunc{p}
      \right)
      \eqnlabel{gausspoly-differentiation} \\
\opname{Integration}
  \noalign{$\sqrt\pi\cdot\int_{-\infty}^{T}\varphi((y,a,b,c),p)(t) \dif t$}
   &=& s\cdot\sqrt\pi\cdot\int_{-\infty}^{T}f(y,a,b,c)(t) \dif t +
       f(y,a,b,c)(T)\cdot \polyfunc{q}(T) \\
   &=& s\cdot\exp\left(-a + \frac{b^2}{4c}\right) \cdot
       \frac{1 + \erf\left(\frac{b}{2\sqrt{c}} + \sqrt{c\pi}\cdot T\right)}
            {2\sqrt{c}} \\
   && \hfill + f(y,a,b,c)(T)\cdot \polyfunc{q}(T) \\
   && \qquad\mbox{where}\quad \polyfunc{q}(t) =
          \frac{\frac{1}{\sqrt\pi}\cdot\polyfunc{q}'(t) - \polyfunc{p}(t) - s}
               {b + c\cdot \sqrt\pi\cdot t}
      \eqnlabel{gausspoly-integration} \\
\eqnbreak
\opname{Eigenfunction of \Fourier{} transform}
   e_n &=& 
       f(1,0,0,2)^{(n)} \cdot f(1,0,0,-1)
      \eqnlabel{gausspoly-eigen}
\end{IEEEeqnarray*}
The equation \eqnref{gausspoly-integration}
is the inverse of \eqnref{gausspoly-differentiation}.
This implies that in \eqnref{gausspoly-integration}
the polynomial~$q$ depends recursively on itself.
However because the degree of $p$ is one more than that of $q$,
the leading term of $q$ only depends on the leading term of $p$.
Thus we can successively determine the terms of $q$ starting at the highest one.
The equation can be translated almost literally
to a polynomial division with remainder $s$
in our Haskell implementation
and be solved using \keyword{lazy evaluation}.

In \eqnref{gausspoly-eigen} we have used
the definition of \person{Hermite} polynomials and the known fact,
that the \Gaussian{} function multiplied with a \person{Hermite} polynomial
is an eigenfunction of the \Fourier{} transform.

\subsection{Mixed \Gaussian{}s}
\seclabel{gaussians4}

In order to support sums of signals,
we must maintain a set $A$ of parameters for the \Gaussian{}s
and a map $P$ from the \Gaussian{} parameters $\alpha$
to the associated polynomial factor.
\[
\sum_{\alpha\in A} f\alpha \cdot \polyfunc{P\alpha}
\]
Eventually, this representation is general enough
in order to be target of a windowed \Fourier{} transform
or a best basis or matching pursuit.
That is we can approximate real world signals in a natural way
and perform exact signal processing operations on them.

\section{Related work}
\seclabel{related-work}

With our paper we wanted to draw a connection
between Computer Algebra on the one side
and Signal Processing and Stochastics on the other side.
With ``Computer Algebra'' we mean exact computations
involving complex mathematical objects
like polynomials, polynomial ideals, groups,
that is at a higher level than computing with individual numbers
but at a lower level than computing with general mathematical expressions
as in symbolic manipulations.
Although signal processing is certainly
not the most prominent application of computer algebra,
there are many problems that were solved using computer algebra methods.
\cite{grabmeier2003computeralgebra}
The most famous application is probably the Discrete \Fourier{} Transform,
that can be performed in log-linear time with respect to the length of the input data
using techniques from number theory, finite fields, polynomial rings
and automated code generation.
\cite{clausen1993fft,frigo2005fftw}
Closely related is the fast convolution of discrete signals,
that uses the Fast \Fourier{} Transform
and by the chirp transform
it is also possible to express a fast \Fourier{} Transform
in terms of fast convolution algorithms.
Another computer algebra application in signal processing
is the design of frequency filters,
where we have to construct rational functions
given conditions for the location of its zeros and poles.
\cite{hamming1989digitalfilter}

In \cite{dalotto1998signalalgebra} the authors develop signal processing
the algebraic way, as we do in our paper.
That is,
opposed to the sample value focus of most signal processing literature
they treat signals as objects,
define operations on them and propose and prove laws.
The book is concerned with two-dimensional signals,
but the difference to one-dimensional signals is not essential in this approach.
Unlike our consideration of real functions
and especially modified \Gaussian{} functions
the authors stick to discrete signals.

Compared to established computer algebra systems
and their symbolic integration machineries,
our framework provides no new class of closedly integrable functions.
All of Maple 10, Mathematica 5.2, Maxima 5.20.1
can integrate the integrals occurring
in convolution, \Fourier{} transform, norm and scalar product
of our products of Gaussian function and polynomial in closed form.
MuPAD 4.0.6 cannot integrate the more complicated convolutions
and Axiom 20091101 cannot cope with the integrals at all,
since it does not yet support assumptions.
We need assertion for the coefficient $c$ of the quadratic term
in the exponent of the \Gaussian{} function.
It must have positive real part,
or must be at least a positive real number.
That is, with an assertion we must exclude signals
with constant amplitude or even unbounded amplitude.



\section{Outlook}
\seclabel{outlook}

\subsection{\person{Dirac} impulses}

The \person{Dirac} impulse $\delta$ is a virtual function
that is infinitely high at the origin and zero elsewhere,
enclosing an infinitely high and infinitely narrow rectangle of area~1.
If this function would exist, it would be the identity element of convolution.
Its \Fourier{} transform is the function that is constant~1,
because this is the identity element of pointwise function multiplication.
In stochastics the \person{Dirac} impulse is needed
for representing mixed discrete/continuous probability distributions.
Existence of derivatives of the \person{Dirac} impulse
would eliminate the need for a distinct differentiation operation,
because $x' = x * \delta'$.
We could also more easily hide the factor $\sqrt{\pi}$
in the differentiation operation
and we could represent frequency spectra of polynomial functions.

Several approaches like \person{Schwartz} distributions
and non-standard analysis were developed,
in order to get a strictly founded notion of a \person{Dirac} impulse.
However, none of them is completely satisfying:
\person{Schwartz} distributions have no longer a notion of function application
and they cannot be multiplied pointwise.
Non-standard analysis allows to define
infinitely high and infinitely narrow functions,
that actually let real functions unaltered
at the ``coarse'' scale of real numbers after convolution.
However when considering a non-standard function on all scales,
convolution with the non-standard \person{Dirac} impulse
well changes the convolution partner.
It is not known to us,
whether an approach can exist at all, that fulfils all expectations.

All the more it is interesting whether we can have an object,
that exactly behaves like a \person{Dirac} impulse in our theory.
Since our theory is abstracted from, but not bound to real functions,
we could check this way, whether a \person{Dirac} impulse makes sense at all.
Formally in our approach a \person{Dirac} impulse
could be represented by $f(1,0,0,+\infty)$.
The term $+\infty$ could be made precise by using projective geometry,
i.e. by allowing an object like $\frac{1}{0}$.
But then it is open, whether we should use
\begin{enumerate}
\item $\exp\left(-\pi\cdot\frac{a+b\cdot t+c\cdot t^2}{d}\right)$,
\item $\exp\left(-\pi\cdot(a+b\cdot t+\frac{c}{d}\cdot t^2)\right)$,
\item $\exp\left(-\pi\cdot(a+b\cdot \frac{t}{d}+c\cdot \frac{t^2}{d^2})\right)$ or
\item $\exp\left(-\pi\cdot\frac{a_0}{a_1}+\frac{b_0}{b_1}\cdot t
         + \frac{c_0}{c_1}\cdot t^2\right)$
\end{enumerate}
with a projective interpretation of the fractions,
and how to cope with the amplitude parameter.

\subsection{Discrete signal processing}

We would like to have the same set of operations and laws
for discrete signals that we already have for real signals.
In order to have dual time and frequency domains,
we need to content ourselves with periodic discrete signals.
For instance for discrete periodic signals $x$ and $y$ of period length $n$
convolution and Discrete \Fourier{} Transform ($\DiscreteFourierTransform$)
are usually defined as:
\begin{IEEEeqnarray*}{rCl}
(x*y)_k &=& \sum_{j\in\Z_n} x_j \cdot y_{k-j} \\
(\DiscreteFourierSynthesis x)_k &=& \frac{1}{\sqrt{n}}\cdot
 \sum_{j\in\Z_n} \exp\left(\frac{2\pi i}{n}\cdot j\cdot k\right) \cdot x_j \\
(\DiscreteFourierAnalysis x)_k &=& \frac{1}{\sqrt{n}}\cdot
 \sum_{j\in\Z_n} \exp\left(-\frac{2\pi i}{n}\cdot j\cdot k\right) \cdot x_j
.
\end{IEEEeqnarray*}
The factor $\frac{1}{\sqrt{n}}$ is chosen,
such that $\DiscreteFourierSynthesis$ becomes unitary.
However with this definition it does not hold
$\DiscreteFourierSynthesis(x\cdot y)
 = 
\DiscreteFourierSynthesis x * \DiscreteFourierSynthesis y$,
but instead
$\sqrt{n}\cdot\DiscreteFourierSynthesis(x\cdot y)
 = 
\DiscreteFourierSynthesis x * \DiscreteFourierSynthesis y$.
One solution would be to add the factor $\frac{1}{\sqrt{n}}$
to the definition of the convolution.
This is at least very uncommon.
An alternative is to treat discrete signals
as piecewise constant functions.
The sums are turned to integrals and thus need a step width.
To this end we equip every signal with a \keyword{sampling rate}
and denote it with $\samplingrate$.
It holds $\samplingrate (\DiscreteFourierSynthesis x) = \frac{n}{\samplingrate x}$.
We obtain the definitions
\begin{IEEEeqnarray*}{rCl}
(x*y)_k &=& \frac{1}{\samplingrate x}\cdot \sum_{j\in\Z_n} x_j \cdot y_{k-j}
 \quad\text{for }\{\samplingrate x, \samplingrate y\} = \{\samplingrate (x*y)\} \\
(\DiscreteFourierSynthesis x)_k
 &=& \frac{1}{\samplingrate x}\cdot
   \sum_{j\in\Z_n} \exp\left(\frac{2\pi i}{n}\cdot j\cdot k\right) \cdot x_j \\
(\DiscreteFourierAnalysis x)_k
 &=& \frac{1}{\samplingrate x}\cdot
   \sum_{j\in\Z_n} \exp\left(-\frac{2\pi i}{n}\cdot j\cdot k\right) \cdot x_j
.
\end{IEEEeqnarray*}

Following the reasoning of real signals,
we need eigenfunctions of the \Fourier{} transform.
Generic simply representable eigenbases of the Discrete \Fourier{} transform
are currently not known,
but according to the \person{Poisson} summation formula,
the discretised and periodically summed eigenfunctions
of the Continuous \Fourier{} transform,
are eigenvectors of the discrete transform.
However discretising an eigenbasis of the continuous \Fourier{} transform
may not yield a discrete eigenbasis.

In discrete signal processing the identity element of convolution
is simple to get:
It is the signal that is 1 at index 0 and zero elsewhere.
In contrast to that, the operation of signal dilation
may lead to undefined and
multiple times defined elements in the resulting vector.
The natural solution is to set undefined elements to zero
and sum up all candidates for multiply defined output elements.
This can be written generally in the following way,
where $n$ is the signal period length
and the empty sum is zero:
\[
(\dilate{k}{x})_j =
   \sum_{l\in\Z_n:\ l\cdot k = j} x_l
.
\]
This definition matches shrinking the vector in the frequency domain.
Nonetheless, we have to drop invertibility of dilation from the list of laws,
that hold for real signals.

Another problem is the definition of differentiation.
We could replace it by centred discrete differences.
Via the \Fourier{} transform this would also yield a notion
of periodic polynomials,
namely polynomials in $\sin\frac{2\pi\cdot t}{n}$ instead of $t$.
But this interpretation of differentiation is different
from differentiation of a continuous function
with subsequent discretisation and periodic summation.
Thus it cannot be used for eigenvector computation
in the same way we used it for continuous signals.


Summarised, we cannot simply perform the operations
on our parameter tuples,
that we developed for continuous signals,
and use them for discrete signals
by just interpreting them in a discrete way.

\bibliographystyle{splncs03}
\bibliography{literature,haskell,audio}

\appendix


\end{document}